\documentclass[a4paper]{article}
\title{New strongly regular graphs derived from the $G_2(4)$ graph}
\author{Thomas Jenrich}
\date{2018-05-09}
\begin{document}
\maketitle

\section{Abstract and introduction}

We consider simple loopless finite undirected graphs. Such a graph is
called strongly regular with parameter set $(v,k,\lambda,\mu)$, or short a
srg$(v,k,\lambda,\mu)$, iff it has exactly $v$ vertices, each of them has
exactly $k$ neighbours, and the number of common neighbours of any two
different vertices is $\lambda$ if they are neighbours and $\mu$
otherwise.

The $G_2(4)$ graph is a well-known srg(416,100,36,20). In this article, we
explicitly construct it and a certain subgraph $E$ induced by 320 vertices in
the same way as in \cite{Jen}. We discover some interesting properties
of $E$ and derive five strongly regular graphs from it:
A srg(256,60,20,12) $F$ which is a subgraph induced by 256 vertices
and four srg(336,80,28,16) $H$, $H_1$, $H_2$ and $H_3$ which do have $E$ as
induced subgraph.

The latter three graphs are new in version 4 of this article and seem to have
been completely unknown as $H$ was before version 1 appeared.
The graph $F$ is isomorphic to objects described in \cite{Bro1}, section
\emph{Subgraphs}, subsection \emph{c) Unions of 16 16-cocliques};
but the strong regularity has been unnoticed before version 1 of this article.

Several propositions in this article have been checked by executing the
additionally (in the source package) provided program G24DGS2 and the program
Dreadnaut from the popular graph theoretic software nauty (by Brendan McKay
and Adolfo Piperno).

\section{Notation}

Universal operator symbols (e.g. the equal sign) are sometimes also used to
connect logical (Boolean) operands.

One can and we will understand a simple graph $G$ on a vertex set $V$ as
a subset of $V \times V$ and express the fact that $x,y \in V$ are adjacent
(are neighbours, are connected (by an edge)) with respect to $G$ by $G(x,y)$.

Using this notation, $G$ is a srg$(v,k,\lambda,\mu)$ iff
 $$|V|=v \land \forall i \in V : \lnot G(i,i) \land |\{j \in V : G(i,j)\}|=k$$

and for all different $i,j \in V$

$$G(i,j)=G(j,i) \land
 |\{ p \in V : G(p,i) \land G(p,j)\}| = \left\{\begin{array}{ll}
\lambda & \mbox{if $G(i,j)$}\\
\mu & \mbox{otherwise}\end{array}\right.$$

For integer numbers $a$ and $b$, $[a,b]$ means the set of integer
numbers from $a$ up to $b$. This set is empty if $a$ exceeds $b$.

\section{A construction of the $G_2(4)$ graph $G$}

The following algorithm is an extract of the description in \cite{Bro1},
refering to \cite{CrnMik}. See also \cite{BrHoKiNa} and \cite{Soi}.

``Consider the projective plane PG(2,16) provided with a nondegenerate
Hermitean form. It has 273 points, 65 isotropic and 208 nonisotropic.
There are 416 = 208 $\cdot$ 12 $\cdot$ 1/6  orthogonal bases. These are the
vertices of $G$. [...] Associated with a basis \{a,b,c\} is the triangle
consisting of the 15 isotropic points on the three lines $ab$, $ac$, and $bc$.
[...] $G$ can be described as the graph on the 416 triangles, adjacent
when they have 3 points in common.''

The construction given above assigns a set of 15 different isotropic points
to each of the 416 vertices/triangles. To be more detailed, it assigns
a set of 5 different isotropic points to each non-isotropic point and 3
different non-isotropic points with pairwise disjoint assigned sets of
isotropic points to each vertex.
For symmetry reasons, from $208 \times 5 = 1040 = 65 \times 16$,
$416 \times 15 = 6240 = 65 \times 96$, and $416 \times 3 = 208 \times 6$,
we conclude that each of the 65 isotropic points is assigned to exactly
16 non-isotropic points and to 96 vertices/triangles, and each of the
208 non-isotropic points is assigned to exactly 6 vertices/triangles.

Assume that integer numbers from 1 to 65 are given to the isotropic points
and consider the sets of the numbers of the isotropic points, here shortly
named iso-sets.
Thus two vertices of $G$ are adjacent if the cardinality of the
intersection of their iso-sets is 3.

\section{Vertex subsets $B$ and $C$ and the subgraph $E$ on 320 vertices}

We define $B$ to be the set of all vertices in $V$ containing the
index 1 in their iso-set and $C$ to be the set of the remaining vertices.
From the previous section, we get $|B|=96$, implying $|C|=416-96=320$.
Let $E$ be the subgraph of $G$ induced by $C$. Consequently,
$E(i,j) = G(i,j)$ for all $i,j \in C$.

The computer program G24DGS2, provided in the source package of this article
and described in more detail in the last sections before the references,
implements the construction of the iso-sets of the members of $V$ and its
just defined subsets, checks also relevant cardinalities.

\section{A four-level hierarchy on $C$ with respect to $E$}

Let $X$ be the set of the 16 non-isotropic points to which the isotropic
point with index 1 is assigned. Consequently, to each vertex in $B$ at
least one element of $X$ is assigned.

For each vertex $i \in C$, let $W(i)$ be the set of those elements of $X$
that are not assigned to any neighbour of $i$ in $B$.

The program G24DGS2 checks that
$$\forall i \in C : |W(i)| = 4\eqno(1)$$
and constructs a bijective function $M$ that for each $r \in [0,4]$ and
$s,t,u \in [0,3]$ returns a vertex $M(r,s,t,u) \in C$
such that for all $r_1,r_2,r_3 \in [0,4]$ and
$s_1,s_2,s_3,t_1,t_2,u_1,u_2 \in [0,3]$ the following propositions hold:

$$
|W(M(r_1,s_1,t_1,u_1)) \cap W(M(r_2,s_2,t_2,u_2))| =
\left\{\begin{array}{rl}
 4 & \mbox{if $r_1 = r_2 \land s_1 = s_2$}\\
 0 & \mbox{if $r_1 = r_2 \land s_1 \ne s_2$}\\
 1 & \mbox{if $r_1 \ne r_2$}\end{array}\right.\eqno(2)
$$

$$ \neg E(M(r_1,s_1,t_1,u_1),M(r_1,s_1,t_1,u_2))\eqno(3)$$

$$ t_1 \ne t_2 \Rightarrow E(M(r_1,s_1,t_1,u_1),M(r_1,s_1,t_2,u_2))\eqno(4)$$

$$ s_1 \ne s_2 \Rightarrow \neg E(M(r_1,s_1,t_1,u_1),M(r_1,s_2,t_2,u_2))\eqno(5)$$

$$ r_1 \ne r_2 \Rightarrow
 |\{ u : u \in [0,3] \land E(M(r_1,s_1,t_1,u_1),M(r_2,s_2,t_2,u))\}|=1\eqno(6)$$

$$ |\{ (t,u) : t,u \in [0,3] \land
   E(M(r_1,s_1,t_1,u_1),M(r_3,s_3,t,u)) \land
   E(M(r_1,s_1,t_2,u_2),M(r_3,s_3,t,u)) \}| $$
$$=
\left\{\begin{array}{rl}
 0 & \mbox{if $t_1=t_2 \land u_1 \ne u_2$}\\
 1 & \mbox{if $t_1 \ne t_2$}\end{array}\right.\eqno(7)
$$

\section{Some reformulations and conclusions}

Here and in the program G24DGS2, we will call

 $\{ M(r,s,t,u) : t,u\in [0,3] \}$, where $r \in [0,4] \land s \in [0,3]$, a cell,

 $\{ M(r,s,t,u) : u \in [0,3] \}$, where $r \in [0,4] \land s,t \in [0,3]$, a cell part,

 $\{ \{ M(r,s,t,u) : t,u\in [0,3] \} : s \in [0,3]\}$, where $r \in [0,4]$, a cell set.

Equations (3) and (4) imply in particular that the subgraph of $E$ induced
by one the 20 cells is a complete fourpartite graph with part size four,
usually denoted by $K_{4,4,4,4}$. And there is no edge between vertices in
different cells within the same cell set.

In each cell we can take a vertex from each part and get a 4-clique. There
are $4^4=256$ such combinations in each cell, 1024 in each cell set.
By (6), each vertex has exactly one neighbour in each of the four parts of a
cell in a different cell set. Together, those five vertices induce a 5-clique.

Each cell part is a 4-coclique. For each cell set, we can take a part
from each of the four cells and those 16 vertices are pairwise non-adjacent.
There are $4^4=256$ such 16-cocliques in a single cell set, 1280 of them
in all five.

For each cell set, we can construct $(4 \times 3 \times 2)^3=13824$
different divisions into four 16-cocliques (by fixing the selection of the
parts in the first cell and combining all permutations of the parts in each
of the three other cells).
For all five sets, that gives $13824^5$ different divisions of $C$ into
twenty 16-cocliques.

Each vertex has exactly $3 \times 4 = 12$ neighbours in its own cell, none in
other cells of the same cell set, and, by (6), exactly 4 neighbours in each
cell in the other four cell sets, summing up to $12+4 \times 4\times 4=76$.

Two different vertices in the same cell part do have $3 \times 4 = 12$ common
neighbours within that cell and, by (5) and (7), no other ones.
Two vertices in different parts of a cell do have $2 \times 4 = 8$ common
neighbours within that cell and, by (5), none in the other cells in the same
set, and, by (7), exactly one common neighbour in each cell in the other four
cell sets, summing up to $8 + 4 \times 4 = 24$.

\section{The induced subgraph $F$ of $E$ on 256 vertices}

We choose four of the five cell sets and consider the subgraph of $E$
induced by the contained 256 vertices. Because there are five such choices,
we get actually five graphs. And the program G24DGS2 does check each of them.
But it also delivered input data for the program Dreadnaut (part of the
software nauty) which then has been used to check that those five graphs are
(pairwise) isomorphic. So we can speak of just one graph and name it $F$.

\subsection{Some properties of $F$}

The only change compared to $E$ is the exclusion of the vertices in one cell
set and of the incident edges. So we have to modify the calculations in the
corresponding subsection just a little bit.

For instance, we can construct $4 \times 256=1024$ different 16-cocliques and
use them to enumerate $13824^4$ different divisions of the vertex set of $F$
into sixteen 16-cocliques.

Each vertex has exactly $12 + 3 \times 4\times 4=60$ neighbours.

Just as for $E$, two different vertices in the same cell part have
12 common neighbours.
Two vertices in different parts of a cell have $8 + 3 \times 4 = 20$.

In order to complete the proof that $F$ is a srg(256,60,20,12), we would have
to estimate the numbers of common neighbours of vertices in different cells
too. An older unpublished version of the program G24DGS2 did just that. But
those checks were removed because they were not as simple as and not
faster than the universal srg check routine that is used now.

\section{The strongly regular supergraph $H$ of $E$ on 336 vertices}

Let $D = \{d(x) : x \in X\}$ be a set of 16 additional vertices not
contained in $V$. Let $H$ be the graph on $C \cup D$ (320+16=336
elements) satisfying
$$\forall i_1,i_2 \in D : \neg H(i_1,i_2)\eqno(8)$$
$$\forall i_1,i_2 \in C : H(i_1,i_2) \Leftrightarrow E(i_1,i_2)\eqno(9)$$
$$\forall i \in C, x \in X : H(i,d(x)) \Leftrightarrow H(d(x),i)
 \Leftrightarrow x \in W(i)\eqno(10)$$

By (8) and another result in a previous section, each vertex $i \in C$ has
exactly 76 neighbours in $C$. By (9), (10), and (1), the only additional
neighbours are four vertices in $D$. Thus, $i$ has exactly 80 neighbours.
For the complete check that $H$ is a srg(336,80,28,16) we refer again to
the dedicated routine call in the program G24DGS2.

\section{The s. r. supergraphs $H_1$, $H_2$ and $H_3$ of $E$ on 336 vertices}

The construction differs from that of $H$ just by permutations
during the use of $W$ for connecting the vertices in C with the
vertices in $D$: The vertices in cell M[4,i] get the connections that are
originally dedicated to the vertices in M[4,p(i)], where p is a permutation
of 0 to 3:

For $H_1$: $p=(0,1,3,2)$.

For $H_2$: $p=(0,2,3,1)$.

For $H_3$: $p=(1,0,3,2)$.

Other permutations (even for other cell sets) work as well,
but do not result in additional non-isomorphic graphs.

The considered graphs on 336 vertices do have different numbers of 6-cliques;
here is an explanation:
The number of 6-cliques in each of the four graphs is the sum of the numbers
of 5-cliques in the neighbourhoods of the 16 vertices in D.
Those neighbourhoods are subgraphs induced by unions of 5 cells from
different cell sets, which contain (as has been checked computationally)
exactly 576 or 640 5-cliques. In fact, when (4) cells from 4 cell sets are
already selected, exactly one cell in the remaining cell set gives the
number 576. The numbers of 6-cliques in the whole graphs are:

$H$: $9024= 16 \times 576$.

$H_1$: $9632= 8 \times 576 + 8 \times 640$.

$H_2$: $9936= 4 \times576 + 12 \times 640$.

$H_3$: $10240= 16 \times 640$.

\section{The provided program G24DGS2}

\subsection{Mathematical foundations of the implementation}

The program uses elements of the three-dimensional vector
space over the finite field GF(16) to constitute and represent the objects
in the projective plane PG(2,16), a very common way.
For an introduction to projective planes, see e.g. \cite{Che}.

The applied Hermitean form takes three-dimensional vectors $a$ and $b$ over GF(16)
and returns \mbox{$a_1 \overline {b_3} + a_2 \overline {b_2} + a_3 \overline {b_1}$},
where addition, multiplication, and conjugation operate in GF(16).

\subsection{Compiling and executing}

The source code file G24DGS2.PAS has been developed for PASCAL compilers
compatible with Turbo Pascal 4.0. Lines are at most 78 characters long. For
inspections the use of an ASCII compatible monospaced font is strongly
recommended. The intended indentation is by one character per structure
level, using blanks (instead of tabs).

The program includes a good portion of comments (enclosed in curly braces).
So it should be fairly understandable at least by readers knowing at least
one imperative programming language.

In order to allow to exclude unwanted (e.g. long running) tasks from
execution without modifying the source text, certain parts of the program
are compiled (and executed) only under the condition that a certain symbol
has been defined during the compilation:

The symbol CHKSRG enables checking the SRG properties of the constructed graphs
that are claimed to be SRGs ($G$, $H$, and all 5 five represenataions
of $F$).

The symbol WRIDRE enables writing of input files for Dreadnaut, for each of
the constructed graphs. Caution: In that case, the program would write
files named G24.DRE, 320.DRE, 336.DRE, 336\_1.DRE, 336\_2.DRE, 336\_3.DRE,
256\_0.DRE, 256\_1.DRE, 256\_2.DRE, 256\_3.DRE, and 256\_4.DRE into the
current (working) directory without explicite confirmation.

The symbol CNTCLI enables the counting of cliques of sizes 2 to 7 for
each of the constructed graphs.

\vspace{0.1in}

The program has been successfully compiled and executed on a 933 MHz Intel
PIII PC running MS Windows 98 SE. These are the used compilers and the
respective four execution times in seconds (roughly measured with some
overhead utilizing the PC clock (resolution: 0.055 s)) from compilations
with none of the three symbols, WRIDRE, CHKSRG, or CNTCLI defined:

Turbo Pascal 5.5 : 3.35 / 4.23 / 58.71 / 152.97

Turbo Pascal 7.01 : 1.37 / 2.20 / 56.52 / 149.29

Borland Delphi 4.0 build 5.37 : 0.77 / 0.99 / 5.11 / 14.23

Virtual Pascal 2.1 build 279 : 0.93 / 1.26 / 9.17 / 23.45

Free Pascal 2.4.4 i386-Win32 : 0.82 / 0.99 / 8.46 / 21.37

\vspace{0.1in}

To avoid a compilation result depending on the settings you could
use the command line versions of the compilers (TPC for Turbo Pascal, BPC
for Borland Pascal 7, DCC32 for Borland Delphi (32 bit versions; do not miss
to use the -CC option in order to generate a console executable), VPC for
Virtual Pascal, FPC for Free Pascal) instead of the compilers integrated in
the IDEs.

The appropriate string to define a symbol for a command line compiler
is usually the concatenation of a certain prefix and the name of that symbol.
The appropriate prefix is -D (alternatively -d) for compilers from Borland
(Turbo Pascal, Delphi), and -d for Free Pascal.
For example, an invokation of Free Pascal could look like

\texttt{fpc G24DGS2 -dCHKSRG -dWRIDRE -dCNTCLI}

In general, the Pascal compilers are not case sensitive with respect to
symbols (and file names).

\subsection{Input and output}

The program ignores any command line parameters or inputs other than
pressing Ctrl-C to cancel the execution.

Despite the optional writing of input files for Dreadnaut mentioned above,
it writes only to the standard output device. In the default case that
will be the monitor screen. But you can redirect the output to a file.

The success of all performed checks and other operations is indicated
by this finishing line:

\texttt{== Regular program stop ==}

These are the displayed clique counts of sizes 2 to 7:

$G$: \texttt{ 2:20800 3:249600 4:873600 5:698880 6:0 7:0}

$E$: \texttt{ 2:12160 3:107520 4:261120 5:129024 6:0 7:0}

$F$: \texttt{ 2:7680 3:51200 4:81920 5:15360 6:0 7:0}

$H$: \texttt{ 2:13440 3:125440 4:330240 5:201984 6:9024 7:0}

$H_1$: \texttt{ 2:13440 3:125440 4:330240 5:201984 6:9632 7:0}

$H_2$: \texttt{ 2:13440 3:125440 4:330240 5:201984 6:9936 7:0}

$H_3$: \texttt{ 2:13440 3:125440 4:330240 5:201984 6:10240 7:0}

\section{Results from Dreadnaut runs}

The binary executed under MS Windows 98 SE on a 933 MHz Intel PIII was
dreadnautB.exe, (according to the starting line) compiled from version
2.2 of Dreadnaut for 32-bit processors, variant BIG, and contained
in the GAP package GRAPE 4r6p1 (downloaded via \cite{GRP}).

In separate runs, each of the 11 files written by G24DGS2 has been
used as input stream (by redirection). In the case of the five graphs on 256
vertices, the full data of the canonically labelled graph has been written
to another (text) file. By simple file comparison, the equality of the
canonically labelled graphs and that way indirectly the isomorphicity of the
five original graphs has been checked.

The following short list gives for each of the graphs a small extract of the
runs: the size of the automorphism group and (enclosed in brackets) the
check sum of the canonically labelled graph :

$G$: \texttt{ 503193600 [4ef1998 7b631ca b27cc78]}

$E$: \texttt{ 368640 [9ef74b4 67d5649 fd37dc0]}

$F$: \texttt{ 368640 [5422323 c601302 4629513]}

$H$: \texttt{ 3840 [6940cf2 40d3533 c70f0b4]}

$H_1$: \texttt{ 128 [2ee14b7 906958b 78981c6]}

$H_2$: \texttt{ 96 [2a0400f 54a4723 38e5992]}

$H_3$: \texttt{ 1280 [de177f3 aab5328 aee8ca1]}

\vspace{0.1in}

Author's eMail address: thomas.jenrich@gmx.de

\end{document}